\documentclass[11pt,reqno]{amsart}
\usepackage[a4paper,margin=27mm]{geometry}
\usepackage{cmap}
\usepackage[OT1]{fontenc}
\usepackage[utf8]{inputenc}
\usepackage[english]{babel}
\usepackage{amsmath,amsfonts,amssymb,mathrsfs,amscd}
\usepackage{mathtools}
\usepackage{textcase}
\usepackage[expansion=false]{microtype}
\usepackage{booktabs}

\usepackage[unicode]{hyperref}
\hypersetup{
  colorlinks=true,
  linkcolor=blue,
  citecolor=blue,
  urlcolor=blue
}

\allowdisplaybreaks[2]
\numberwithin{equation}{section}

\theoremstyle{plain}
\newtheorem{theorem}{Theorem}[section]
\newtheorem{lemma}[theorem]{Lemma}

\newtheorem{corollary}{Corollary}[section]

\date{}

\begin{document}
\title{The Estermann Problem
with Almost Proportional Summands}
\author{K.~I.~Mirzoabdughafurov}
\address{Tajik State University of Finance and Economics,
Dushanbe, Tajikistan}
\email{karimjon2003@mail.ru}

\begin{abstract}
We obtain a family of estimates for short Weyl sums
$$
\left|
\sum_{x-y<m\le x}e^{2\pi i\alpha m^n}
\right|
\ll
y\left(
\frac{1}{q}+\frac{1}{y}+\frac{q}{y^n}
\right)^{\frac{s-1}{s\cdot2^{n-1}}}
\bigl(\log(qy)\bigr)^{
\frac{(n-1)^s+s-2}{s\cdot2^{n-1}}},
$$
where $n\ge3$ and $s\ge2$ are fixed integers,
$x\ge y\ge2$, $\alpha\in\mathbb R$,
$a\in\mathbb Z$, $q\in\mathbb N$,
$(a,q)=1$, and
$\left|\alpha-\frac{a}{q}\right|\le\frac{1}{q^2}$.

For $n\ge4$, we also derive from the modern form
of Vinogradov's mean value theorem the estimate
$$
\left|
\sum_{x-y<m\le x}e^{2\pi i\alpha m^n}
\right|
\ll
y\left(
\frac{1}{q}+\frac{1}{y}+\frac{q}{y^n}
\right)^{\frac{1}{n(n-1)+2}}
\log y.
$$

Using these estimates, repeated Dirichlet
approximation, and an estimate in the intermediate range,
we establish an asymptotic formula
for the number of representations
$$
p_1+p_2+m^n=N,\quad
|p_i-\mu_iN|\le H,\quad i=1,2,
\quad |m^n-\mu_3N|\le H,
$$
where $p_1,p_2$ are primes,
$m$ is a positive integer,
and $\mu_1,\mu_2,\mu_3$ are fixed
positive numbers whose sum is $1$.
The asymptotic formula holds in the range
$$
N^{1-\theta_n}(\log N)^{A_n}\le H\le N^{1-\delta},
\qquad
\theta_n=\min\left\{\frac{1}{2n},\frac{2}{n^2}\right\},
$$
where $0<\delta<\theta_n$ is fixed,
$A_3=3$, $A_4=\frac{61}{4}$, and
$A_n=6n-12+\frac{22}{n}-\frac{12}{n^2}$
for $n\ge5$.
In the cubic case, the logarithmic exponent
in the previously known lower bound for $H$
is reduced from $10$ to $3$.
For $n\ge4$, the exponent of $N$
is reduced from $1-\frac{1}{n(n-1)}$
to $1-\frac{2}{n^2}$.
\end{abstract}
\subjclass[2020]{Primary 11P32; Secondary 11L07, 11P55}

\keywords{Estermann problem, almost proportional
summands, short Weyl sums, Vinogradov's mean value
theorem, repeated Dirichlet approximation,
intermediate range of rational approximations,
Hardy--Littlewood method}

\maketitle

\section{Introduction}
\begin{quote}
{\footnotesize
\noindent\textbf{Notation.}
We denote by $\mathbb{N}$, $\mathbb{Z}$,
and $\mathbb{P}$ the sets of
positive integers, integers, and primes, respectively.
Throughout, we put
\[
e(\alpha)=e^{2\pi i\alpha},
\qquad
S(a,q)=\sum_{k=1}^{q}
e\left(\frac{ak^n}{q}\right),
\qquad
a\in\mathbb Z,\quad q\in\mathbb N,
\]
\[
T(\alpha;x,y)
=
\sum_{x-y<m\le x}e(\alpha m^n),\quad
S_1(\alpha;z,h)
=\sum_{z-h<n\le z}\Lambda(n)e(\alpha n),
\quad 0<h\le z,
\]

where $\Lambda(n)$ is the von Mangoldt function.

We write $\|t\|$ for the distance from $t$ to the nearest
integer and $[t]$ for the integer part of $t$.

The functions $\mu(q)$ and $\varphi(q)$ are the M\"obius
and Euler functions, respectively.

The symbol $\varepsilon$ denotes a sufficiently small
fixed positive number.}
\end{quote}
Let $n\ge3$ be a fixed integer,
and let $\mu_1,\mu_2,\mu_3$ be fixed positive
numbers satisfying
\begin{equation}\label{eq:proportions}
\mu_1+\mu_2+\mu_3=1.
\end{equation}
For a positive integer $N$ and a real number $H>0$,
let $J(N,H)$ denote the number of ordered triples
$(p_1,p_2,m)$ satisfying the equation
\begin{equation}\label{eq:equation}
p_1+p_2+m^n=N,
\end{equation}
where $p_1,p_2$ are primes and $m$ is a positive
integer, subject to
\begin{equation}\label{eq:windows}
|p_k-\mu_kN|\le H,\qquad k=1,2,
\qquad
|m^n-\mu_3N|\le H.
\end{equation}
When $\mu_1=\mu_2=\mu_3=\frac{1}{3}$, this becomes
the Estermann problem with almost equal summands.
Throughout, we put $L=\log N$.

In 1937, T.~Estermann~\cite{Estermann} proved an asymptotic
formula for the number of representations of a sufficiently large
positive integer as the sum of two primes
and a square.
The cubic version of this problem with almost equal
summands was studied in
\cite{Cubic2014,RakhmonovFozilova2012}.
In particular, Z.~Kh.~Rakhmonov~\cite{Cubic2014}
obtained an asymptotic formula for
\[
H\ge N^{\frac{5}{6}}L^{10}.
\]
Estimates for short cubic Weyl sums were obtained in \cite{RakhmonovMirzo2008}.

Subsequently, F.~Z.~Rakhmonov~\cite{Estermann2024}
considered equation~\eqref{eq:equation}
with the restrictions~\eqref{eq:windows}
for an arbitrary fixed $n\ge3$
and established an asymptotic formula
with the lower bound
\[
H\ge
N^{1-\frac{1}{n(n-1)}}
L^{\frac{2^{n+1}}{n-1}+n-1}.
\]
For $n=3$, this bound also takes the form
$H\ge N^{\frac{5}{6}}L^{10}$.

A substantial analytic difficulty in problems
of this type concerns estimates for short Weyl sums.
In the problem under consideration,
the restrictions~\eqref{eq:windows} lead to the sum
$T(\alpha;x,y)$ with parameters
\[
x\asymp N^{\frac{1}{n}},
\qquad
y\asymp\frac{H}{N^{\frac{n-1}{n}}}.
\]

For comparison, we recall an estimate
of F.~Z.~Rakhmonov~\cite{Estermann2024}.
Let
\[
x\ge x_0,
\qquad
\log x<y\le\frac{x}{\log x},
\]
\[
\alpha=\frac{a}{q}+\lambda,
\qquad
a\in\mathbb Z,\quad q\in\mathbb N,
\qquad
(a,q)=1,
\qquad
|\lambda|\le\frac{1}{q^2}.
\]
Then
\[
|T(\alpha;x,y)|
\ll
y\left(
\frac{1}{q}+\frac{1}{y}+\frac{q}{y^n}
\right)^{\frac{1}{2^n}}
\bigl(\log(qy)\bigr)^{\frac{(n-1)^2}{2^n}}.
\]
In this paper, we obtain a family of estimates
in which the relation between the power
and logarithmic exponents is controlled by
an integer parameter $s\ge2$.

Our first result is stated as
follows.
\begin{theorem}\label{thm:log-weyl}
Let $n\ge3$ and $s\ge2$ be fixed integers,
$x\ge y\ge2$, $\alpha\in\mathbb R$, $a\in\mathbb Z$,
$q\in\mathbb N$, $(a,q)=1$, and
\[
\alpha=\frac{a}{q}+\lambda,
\qquad
|\lambda|\le\frac{1}{q^2}.
\]
Then the sum
\[
T(\alpha;x,y)=\sum_{x-y<m\le x}e(\alpha m^n)
\]
satisfies the estimate
\begin{equation}\label{eq:log-weyl}
|T(\alpha;x,y)|
\ll
y\left(
\frac{1}{q}+\frac{1}{y}+\frac{q}{y^n}
\right)^{\frac{s-1}{s\cdot2^{n-1}}}
\bigl(\log(qy)\bigr)^{\frac{(n-1)^s+s-2}{s\cdot2^{n-1}}}.
\end{equation}
For $n\ge4$, we also have
\begin{equation}\label{eq:vmvt-log}
|T(\alpha;x,y)|
\ll
y\left(
\frac{1}{q}+\frac{1}{y}+\frac{q}{y^n}
\right)^{\frac{1}{n(n-1)+2}}
\log y.
\end{equation}
The implied constant
in~\eqref{eq:log-weyl} depends only on $n$ and $s$,
and that in~\eqref{eq:vmvt-log} depends only on $n$.
\end{theorem}

The proof of~\eqref{eq:log-weyl}
uses $n-1$ successive
differencing steps.
Grouping the terms according to the product
of the differencing parameters leads to a sum
involving the function $\tau_{n-1}(h)$, which counts
the ordered factorizations of $h$ into $n-1$
positive factors.
Applying H\"older's inequality with exponents
$s$ and $\frac{s}{s-1}$ and an estimate for
the $s$th moment of $\tau_{n-1}(h)$
yields the exponents
in~\eqref{eq:log-weyl}.

For $s=2$, the exponents of the parenthesized factor and the logarithm
in~\eqref{eq:log-weyl} coincide with those
in the estimate from~\cite{Estermann2024} quoted above,
while the restriction on the interval length is replaced
by $x\ge y\ge2$.
For $s>2$, the exponent of the parenthesized factor increases,
as does the exponent of the logarithm.
Thus, the parameter $s$ allows a choice
of the balance between the power saving
and the logarithmic loss according to
the range of parameters.

In particular, for $n=3$ and $s=3$, we obtain
\[
|T(\alpha;x,y)|
\ll
y\left(
\frac{1}{q}+\frac{1}{y}+\frac{q}{y^3}
\right)^{\frac{1}{6}}
\bigl(\log(qy)\bigr)^{\frac{3}{4}}.
\]
In their common range of applicability, this estimate
is stronger than the estimate from~\cite{Estermann2024}
provided that
\[
\left(
\frac{1}{q}+\frac{1}{y}+\frac{q}{y^3}
\right)
\bigl(\log(qy)\bigr)^6
\rightarrow0.
\]

Estimate~\eqref{eq:vmvt-log} follows
from the modern form of Vinogradov's mean value theorem
\cite[Corollary~1.3]{Wooley2019}
by the argument
given in~\cite[p.~1608]{Wooley2012}.
It complements the family~\eqref{eq:log-weyl}
and is used for larger values of $n$.

In the application to the Estermann problem,
for $n=3,4$ we use
estimate~\eqref{eq:log-weyl} with $s=3$,
and for $n\ge5$ we use~\eqref{eq:vmvt-log}.
The choice between these estimates determines
the exponent $B_n$ in $Q=L^{B_n}$
and enters the explicit value of~$A_n$.

In the author's paper~\cite{Waring9Cubes} on
Waring's problem for nine cubes with almost equal
summands, an approach to estimating the minor arc integral
was developed, based on repeated rational
approximation by Dirichlet's theorem and an estimate for a short cubic
sum in the intermediate range.
Further developments of this approach
for fourth powers and arbitrary fixed
powers are presented
in the preprints~\cite{IRFourth,IR}.

In~\cite{IR}, an intermediate range estimate is obtained
for an arbitrary fixed degree $n\ge3$,
and Waring's problem with almost
proportional summands is studied.
In the present paper, this approach is applied
to the Estermann problem together
with the new family of estimates~\eqref{eq:log-weyl}.

The first Dirichlet approximation divides the domain
of integration into major and minor arcs.
To estimate the minor arc integral, Dirichlet's theorem
is applied again with the parameter
\[
P=\frac{y^n}{Q},
\quad Q=L^{B_n},
\]
where the exponent $B_n$ is defined below.
This gives the representation
\[
\alpha=\frac{b}{r}+\beta,
\qquad
(b,r)=1,
\qquad
1\le r\le P,
\qquad
|\beta|\le\frac{1}{rP}.
\]
For $r>Q$, Theorem~\ref{thm:log-weyl}
gives the required logarithmic saving.
For $r\le Q$, the denominator $r$ may be small,
and this estimate alone does not provide
such a saving.
Here we use additional information:
the condition that $\alpha$ lies on the minor arcs
implies a lower bound for $|\beta|$.
Together with the second approximation, this gives
\[
\frac{1}{rx^{n-2}y}
\ll|\beta|
\le\frac{Q}{ry^n}
\le\frac{1}{ry^{n-1}}.
\]
In this intermediate range, we apply
the author's estimate from~\cite[formula~(8)]{IR},
stated below as Lemma~\ref{ShortExposumIR}.
In applying it, we retain
both bounds for $|\beta|$, including the sharper
upper bound $\frac{Q}{ry^n}$.
This yields an estimate for $T(\alpha;x,y)$
even when the denominator of the second approximation is small.

For $n\ge4$, the resulting term
\[
Q^{1-\frac{1}{n}}x^{\frac{n-2}{2}}y^{1-\frac{n}{2}}
\]
is of order at most $yL^{-2}$ provided that
\[
y\ge
x^{1-\frac{2}{n}}
L^{\frac{4}{n}+\frac{2(n-1)B_n}{n^2}}.
\]
This scale for the length of the short sum
corresponds to the lower bound
$H\ge N^{1-\frac{2}{n^2}}L^{A_n}$
in Theorem~\ref{thm:main}.
Thus, controlling the error in the second
approximation determines the power range,
whereas estimates~\eqref{eq:log-weyl}
and~\eqref{eq:vmvt-log} allow us to optimize the choice of the exponent $B_n$
in the parameter $Q=L^{B_n}$.

For the outer part of the major arcs, we use
Corollary~3 of~\cite{RakhmonovWaring2024}.
The estimates for these parts of the integral are combined
in the condition~\eqref{eq:Hrange}.

To state the second result,
let $\rho(N,p)$ denote the number of solutions
of the congruence
\[
u^n\equiv N\pmod p
\]
modulo a prime $p$, and put
\begin{equation}\label{eq:euler-product}
\mathfrak S(N)
=
\prod_p
\left(
1+\frac{\rho(N,p)-1}{(p-1)^2}
\right).
\end{equation}
For each integer $n\ge3$, put
\[
\theta_n=
\min\left\{\frac{1}{2n},\frac{2}{n^2}\right\},
\qquad
B_n=
\begin{cases}
\frac{33}{2}, & n=3,\\[4pt]
38, & n=4,\\[4pt]
3\bigl(n(n-1)+2\bigr), & n\ge5,
\end{cases}
\]
and
\[
A_n=
\begin{cases}
3, & n=3,\\[4pt]
\frac{4}{n}+\frac{2(n-1)B_n}{n^2}, & n\ge4.
\end{cases}
\]
Then $B_n>2n$ and $A_n\ge n$.

\begin{theorem}\label{thm:main}
Let $n\ge3$ be a fixed integer,
let $\mu_1,\mu_2,\mu_3$ be fixed positive numbers
satisfying~\eqref{eq:proportions},
and let $0<\delta<\theta_n$ be fixed.
Then, for every sufficiently large positive integer $N$,
uniformly for real $H$ in the range
\begin{equation}\label{eq:Hrange}
N^{1-\theta_n}L^{A_n}
\le H\le N^{1-\delta},
\qquad L=\log N,
\end{equation}
we have the asymptotic formula
\begin{equation}\label{eq:main}
J(N,H)
=
\frac{3\mathfrak S(N)H^2}
{n\mu_3^{\frac{n-1}{n}}N^{\frac{n-1}{n}}L^2}
+
O\left(
\frac{H^2}{N^{\frac{n-1}{n}}L^3}
\right).
\end{equation}
The implied constant in the error term may depend only
on $n$, $\mu_1,\mu_2,\mu_3$, and $\delta$.
The product~\eqref{eq:euler-product} converges absolutely,
and, uniformly in $N$, we have
\[
0<c\le\mathfrak S(N)\le C_n,
\]
where $c$ is an absolute constant
and $C_n$ depends only on $n$.
\end{theorem}
For $n=3$, the condition~\eqref{eq:Hrange}
takes the form
\[
N^{\frac{5}{6}}L^3\le H\le N^{1-\delta}.
\]
Thus, the logarithmic exponent $10$
in the previously known lower bound for $H$
is replaced by $3$.
For each fixed $n\ge4$,
Theorem~\ref{thm:main} reduces the exponent of $N$
in the lower bound for $H$ from $1-\frac{1}{n(n-1)}$ in~\cite{Estermann2024}
to $1-\frac{2}{n^2}$.

In particular, for fourth powers, the exponent
of $N$ in the lower bound for $H$
is reduced from $\frac{11}{12}$ to $\frac{7}{8}$.
The asymptotic formula holds
for arbitrary fixed positive
proportions of the summands, including the case
of almost equal summands.
\begin{corollary}\label{cor:equal}
Let $n\ge3$ be a fixed integer, let
$\mu_1=\mu_2=\mu_3=\frac{1}{3}$
and suppose that the hypotheses of Theorem~\ref{thm:main} hold.
Then the number $J(N,H)$ of solutions
of equation~\eqref{eq:equation} subject to
\[
\left|p_i-\frac{N}{3}\right|\le H,
\qquad i=1,2,
\qquad
\left|m^n-\frac{N}{3}\right|\le H
\]
satisfies the asymptotic formula
\[
J(N,H)
=
\frac{3^{2-\frac{1}{n}}\mathfrak S(N)H^2}
{nN^{\frac{n-1}{n}}L^2}
+
O\left(
\frac{H^2}{N^{\frac{n-1}{n}}L^3}
\right).
\]
The implied constant in the error term may depend
only on $n$ and $\delta$.
\end{corollary}

\section{Auxiliary results}

\begin{lemma}{\rm({\bf Dirichlet})} {\rm \cite[Lemma~2, p.~158]{KarOATCh}.}\label{lemmaDirichle}
For any real number \(\alpha\) and any parameter \(P\geq 1\),
there exist coprime integers \(a\) and \(q\), $1\leq q\leq P,$
such that
\[
\left|\alpha-\frac{a}{q}\right|\leq \frac{1}{qP}.
\]
\end{lemma}
\begin{lemma}{\rm \cite[Lemma~1, p.~12]{Korobov1989}}\label{linear-sum}
Let $\alpha$ be an arbitrary real number,
let $Q$ be an integer, and let $P$ be a positive integer. Then
\[
\left|
\sum_{x=Q+1}^{Q+P} e^{2\pi i\alpha x}
\right|
\le
\min\left(P,\frac{1}{2\|\alpha\|}\right),
\]
where $\|\alpha\|$ is the distance from $\alpha$ to the nearest
integer.
\end{lemma}

\begin{lemma}{\rm \cite[Lemma~2.3, p.~18]{Von}}\label{Weyl}
For a function $\varphi$, put
\[
\Delta_1(\varphi(x);h)=\varphi(x+h)-\varphi(x)
\]
and, for $j\ge1$, define
\[
\Delta_{j+1}(\varphi(x);h_1,\ldots,h_{j+1})=\Delta_j(\varphi(x+h_{j+1});h_1,\ldots,h_j)-\Delta_j(\varphi(x);h_1,\ldots,h_j).
\]
Let
\[
T(\varphi)=\sum_{x=1}^{Q}e\bigl(\varphi(x)\bigr),
\]
where $\varphi$ is an arbitrary arithmetic function. Then
\[
|T(\varphi)|^{2^j}
\le
(2Q)^{2^j-j-1}
\sum_{|h_1|<Q}\cdots\sum_{|h_j|<Q}T_j,
\]
where
\[
T_j=
\sum_{x\in I_j}
e\bigl(\Delta_j(\varphi(x);h_1,\ldots,h_j)\bigr),
\]
and the possibly empty intervals $I_j=I_j(h_1,\ldots,h_j)$
satisfy
\[
I_1(h_1)\subset[1,Q],
\qquad
I_j(h_1,\ldots,h_j)
\subset
I_{j-1}(h_1,\ldots,h_{j-1}).
\]
\end{lemma}

\begin{lemma}\label{orsumm} {\rm{\cite[Theorem~2, p.~61]{Archipov}}.}
Let $n\ge3$ be an integer, let
\[
f(x)=a_nx^n+\cdots+a_1x+a_0
\]
be a polynomial with integer coefficients,
and let $q$ be a positive integer with
\[
(a_n,\ldots,a_1,q)=1.
\]
Then
\[
\left|
\sum\limits_{x=1}^{q}e\left(\frac{f(x)}{q}\right)
\right|
\le c(n)q^{1-\frac{1}{n}},
\]
where
\[
c(n)=
\begin{cases}
e^{4n}, & n\ge10,\\[2mm]
e^{nA(n)}, & 3\le n\le9,
\end{cases}
\]
and
\[
\begin{aligned}
A(3)&=6{,}1, & A(4)&=5{,}5, & A(5)&=5, & A(6)&=4{,}7,\\
A(7)&=4{,}4, & A(8)&=4{,}2, & A(9)&=4{,}05.
\end{aligned}
\]
\end{lemma}
\begin{lemma}\label{ShortExposumIR}
{\rm \cite[formula~(8)]{IR}.}
Let $n\ge3$ be a fixed integer,
\[
x\ge x_0>0,\qquad 0<y\le0.01x,
\]
and let
\[
\alpha=\frac{a}{q}+\lambda,\qquad
a\in\mathbb Z,\qquad q\in\mathbb N,\qquad
(a,q)=1.
\]
Then
\[
\begin{aligned}
\left|T\left(\alpha;x,y\right)\right|
\ll &
q^{1-\frac{1}{n}}\ln q
+
q^{1-\frac{1}{n}}
\ln\left(2+|\lambda|x^{n-1}\right)
\\
&+
q^{-\frac{1}{n}}
\left(|\lambda|x^{n-2}\right)^{-\frac{1}{2}}
+
q^{1-\frac{1}{n}}y
\left(|\lambda|x^{n-2}\right)^{\frac{1}{2}}.
\end{aligned}
\]
The implied constant depends only on $n$.
\end{lemma}
In Lemmas \ref{exposumm1} and \ref{exposumm2}, we assume that $x\ge x_0$,
$1\le y\le 0{,}01x$, and
\[
\alpha=\frac{a}{q}+\lambda,\qquad
a\in\mathbb Z,\quad q\in\mathbb N,\quad
(a,q)=1,\quad q\le\tau,\quad
|\lambda|\le\frac{1}{q\tau}.
\]
\begin{lemma}{\rm\cite[Corollary~2]{RakhmonovWaring2024}.}\label{exposumm1}\quad {\it Let $\tau\geq 2n (n-1)x^{n-2}y$, $q\leq \tau $,
$\alpha=\frac{a}{q}+\lambda$; $(a,q)=1$, $|\lambda|\le \frac{1}{2nq
x^{n-1}}$. Then we have }
$$
T(\alpha,x,y)=\frac{y}{q}S(a,q)\gamma(\lambda;x,y)+O(q^{\frac{1}{2}+\varepsilon
}),
$$
where
$$
\gamma(\lambda;x,y)=\int\limits_{-\frac{1}{2}}^{\frac{1}{2}}e\Bigl(\lambda(x-\frac{y}{2}+yt)^n\Bigr)dt.
$$
\end{lemma}
\begin{lemma}{\rm\cite[Corollary~3]{RakhmonovWaring2024}.}
\label{exposumm2}
Let
\[
\tau\ge2n(n-1)x^{n-2}y,
\qquad
\frac{1}{2nqx^{n-1}}<|\lambda|\le\frac{1}{q\tau}.
\]
Then we have the estimate
\[
T(\alpha;x,y)
\ll
q^{1-\frac{1}{n}}\ln q
+
\min\left(
yq^{-\frac{1}{n}},
x^{\frac{1}{2}}q^{\frac{1}{2}-\frac{1}{n}}
\right).
\]
\end{lemma}

The implied constant depends only on $n$.

\begin{lemma}{\rm\cite[Theorem~1, p.~612]{Sobirov}.}\label{lem:primes-input}
Let $x\ge x_0$, $\mathcal L=\log(xq)$, and let $A$ and $b$ be
arbitrary fixed positive numbers.
Suppose also that
\[
\alpha=\frac{a}{q}+\lambda,\qquad
a\in\mathbb Z,\quad q\in\mathbb N,\quad
(a,q)=1,\qquad 1\le q\le\mathcal L^b.
\]
Then, for
$|\lambda|\le x(2\pi y^2)^{-1}$ and
$y\ge x^{\frac{5}{8}}\mathcal{L}^{1.5A+0.25b+18}$
we have
\[
S_1(\alpha;x,y)
=
\frac{\mu(q)}{\varphi(q)}
\frac{\sin \pi\lambda y}{\pi\lambda}
e\left(\lambda\left(x-\frac{y}{2}\right)\right)
+
O\left(y\mathcal{L}^{-A}\right).
\]
\end{lemma}
\begin{lemma}{\rm ~\cite[p.~174]{UVatson}.}\label{vatson} For $m>0$ and $n\in\mathbb N$, we have the identity
\begin{align*}
&\int\limits_0^{\infty}\frac{\sin^n(mt)}{t^n}dt=\frac{\pi
m^{n-1}}{2^n(n-1)!}\left[n^{n-1}-\frac{n}{1!}(n-2)^{n-1}+\right.\\
 &\left.+\frac{n(n-1)}{2!}(n-4)^{n-1}
-\frac{n(n-1)(n-2)}{3!}(n-6)^{n-1}+\ldots\right].
\end{align*}
\end{lemma}
\section{Proof of Theorem \ref{thm:log-weyl}}

Put
\[
u_0=[x-y],
\qquad
M=[x]-[x-y],
\qquad
k=n-1.
\]
Then $M\ge2$ and $M\asymp y$.

Making the substitution $m=u_0+v$, we obtain
\[
T(\alpha,x,y)=\sum_{v=1}^{M}e(F(v)),
\]
where
\[
F(v)=\alpha(u_0+v)^n
=\alpha v^n+\beta_1v^{n-1}+\cdots+\beta_n.
\]
Here $\beta_1,\ldots,\beta_n$ are real numbers.
All subsequent estimates are uniform in the lower
coefficients of the polynomial $F$.

We apply Lemma~\ref{Weyl} with
$Q=M$, $\varphi(v)=F(v)$, and $j=k=n-1$.
Since the leading coefficient of the polynomial $F$
is $\alpha$, we have
\[
\Delta_k(F(v);h_1,\ldots,h_k)
=
n!\alpha h_1\cdots h_kv+C(h_1,\ldots,h_k),
\]
where $C(h_1,\ldots,h_k)$ is independent of $v$.
Taking the constant factors of modulus one
outside the inner sums and applying the triangle
inequality, we obtain
\begin{align*}
|T(\alpha,x,y)|^{2^k}
\ll &
M^{2^k-k-1}
\sum_{|h_1|<M}\cdots
\sum_{|h_k|<M}
\left|
\sum_{v\in I_k}
e(n!\alpha h_1\cdots h_kv)
\right|.
\end{align*}

We first separate the tuples in which at least one
of the shifts $h_1,\ldots,h_k$ is zero.
The number of such tuples is
\[
(2M-1)^k-(2M-2)^k\ll  M^{k-1}.
\]
For each of them, the inner sum
is at most $M$.
Hence their total contribution is bounded
by
\[
M^{2^k-k-1}\cdot O (M^{k-1})\cdot M
\ll   M^{2^k-1}.
\]

For the remaining tuples, we apply
Lemma~\ref{linear-sum}.
Since each interval $I_k$
contains at most $M$ integers, we have
\[
\left|
\sum_{v\in I_k}
e(n!\alpha h_1\cdots h_kv)
\right|
\le
\min\left(
M,
\frac{1}{2\|n!\alpha h_1\cdots h_k\|}
\right).
\]
The right-hand side depends only on the absolute values
of the shifts $h_1,\ldots,h_k$.
Putting $r_i=|h_i|$, $1\le i\le k$,
and absorbing the factor $2^k$ into the implied constant,
we obtain
\begin{equation*}
|T(\alpha,x,y)|^{2^k}
\ll
M^{2^k-1}+
M^{2^k-k-1}
\sum_{r_1=1}^{M-1}\cdots
\sum_{r_k=1}^{M-1}
\min\left(
M,
\frac{1}{2\|n!\alpha r_1\cdots r_k\|}
\right).
\end{equation*}

Put
\[
v_h=
\min\left(M,\frac{1}{2\|n!\alpha h\|}\right)
\]
and let
\[
\tau_r(h)
=
\sum_{\substack{u_1,\ldots,u_r\ge1\\
u_1\cdots u_r=h}}1
\]
denote the number of ordered factorizations of $h$
into $r$ positive factors.
The number of tuples $r_1,\ldots,r_k$
with a given product $h$
is at most $\tau_k(h)$.
Thus, grouping the terms according to the product
$h=r_1\cdots r_k$, we obtain
\begin{equation}\label{eq:two-differences}
|T(\alpha,x,y)|^{2^k}
\ll 
M^{2^k-1}
+
M^{2^k-k-1}
\sum_{h\le M^k}\tau_k(h)v_h.
\end{equation}
We estimate the sum $\sum\limits_{h\le M^k}v_h$.
Put
\[
D=\frac{1}{q}+\frac{1}{M}+\frac{q}{M^n},
\qquad
l=\log(2qM).
\]
Write
\[
g=(q,n!),
\qquad
q_1=\frac{q}{g},
\qquad
a_1=\frac{n!a}{g},
\qquad
\theta=n!\lambda.
\]
Then
\[
n!\alpha=\frac{a_1}{q_1}+\theta,
\qquad
(a_1,q_1)=1,
\]
\[
\frac{q}{n!}\le q_1\le q,
\qquad
|\theta|\le\frac{n!}{q^2}\le\frac{n!}{q_1^2}.
\]

We divide the integers $1,\ldots,M^k$ into blocks of
\[
R=\left[\frac{q_1}{2n!}\right]+1
\]
consecutive integers; the last block
may be incomplete.
Let $\mathcal B$ be one of these blocks.
If $h_1,h_2\in\mathcal B$ and $h_1\ne h_2$, then
\[
0<|h_1-h_2|
\le R-1
\le\frac{q_1}{2n!}
<q_1.
\]
The condition $(a_1,q_1)=1$ implies that
\[
\left\|
\frac{a_1(h_1-h_2)}{q_1}
\right\|
\ge\frac{1}{q_1}.
\]
Hence
\begin{align*}
\|n!\alpha(h_1-h_2)\|=
\left\|
\frac{a_1(h_1-h_2)}{q_1}
+\theta(h_1-h_2)
\right\|&\ge
\left\|
\frac{a_1(h_1-h_2)}{q_1}
\right\|
-
|\theta|\,|h_1-h_2|\\
&\ge
\frac{1}{q_1}
-
\frac{n!}{q_1^2}\cdot\frac{q_1}{2n!}
=
\frac{1}{2q_1}.
\end{align*}
Consequently, every half-open interval of length
$\frac{1}{2q_1}$ contains at most one
fractional part of $n!\alpha h$ with $h\in\mathcal B$.

In particular, the condition
\[
\|n!\alpha h\|<\frac{1}{2q_1}
\]
is satisfied by at most two integers
$h\in\mathcal B$.
Their contribution is at most $2M$.

For each $j=1,\ldots,q_1$, the condition
\[
\frac{j}{2q_1}
\le\|n!\alpha h\|
<
\frac{j+1}{2q_1}
\]
is also satisfied by at most two integers
$h\in\mathcal B$.
For these integers,
\[
v_h\le
\frac{1}{2\|n!\alpha h\|}
\le\frac{q_1}{j}.
\]
Thus,
\[
\sum_{h\in\mathcal B}v_h
\le
2M+2q_1\sum_{j=1}^{q_1}\frac{1}{j}
\ll
M+q_1\log(2q_1).
\]

Since $R\ge\frac{q_1}{2n!}$,
the number of blocks is at most
\[
\frac{M^k}{R}+1
\le
\frac{2n!M^k}{q_1}+1.
\]
Summing over the blocks and using $k=n-1$, we obtain
\begin{align}
\sum_{h\le M^k}v_h
&\ll
\left(\frac{M^k}{q_1}+1\right)
\bigl(M+q_1\log(2q_1)\bigr)\ll
\left(\frac{M^k}{q}+1\right)
\bigl(M+q\log(2q)\bigr)\nonumber\\
&\ll 
\frac{M^n}{q}
+M^{n-1}\log(2q)
+M+
q\log(2q)\ll  M^nDl.
\label{eq:spacing}
\end{align}

We now estimate the sum
$\sum\limits_{h\le Z}\tau^s_k(h)$, where $Z\ge2$.
For any positive integers $u,v,r$, we have
\[
\tau_r(uv)\le\tau_r(u)\tau_r(v).
\]
Indeed, let $p$ be a prime,
and let $m,n$ be nonnegative integers.
By definition, $\tau_r(p^{m+n})$ is the number
of ordered factorizations
\[
p^{m+n}=u_1\cdots u_r,
\qquad u_1,\ldots,u_r\in\mathbb N.
\]
Each factor $u_i$ is a power of $p$,
that is, $u_i=p^{a_i}$, where $a_i\ge0$.
Thus, each such factorization
corresponds uniquely to a tuple of nonnegative
integers $a_1,\ldots,a_r$ satisfying
the equation
\[
a_1+\cdots+a_r=m+n.
\]
Hence
\[
\tau_r(p^{m+n})
=
\binom{m+n+r-1}{r-1}.
\]
Expanding the binomial coefficient, we obtain
\begin{equation*}
\tau_r(p^{m+n})
=\frac{(m+n+r-1)!}{(m+n)!(r-1)!}=
\prod_{j=1}^{r-1}\frac{m+n+j}{j}\le\prod_{j=1}^{r-1}
\frac{(m+j)(n+j)}{j^2}
=\tau_r(p^m)\tau_r(p^n),
\end{equation*}
since
\[
(m+j)(n+j)-j(m+n+j)
=m n\ge0.
\]
The inequality for arbitrary $u,v$
follows from the multiplicativity of $\tau_r$.

We prove by induction on $s$ that
\[
\tau^s_k(h)\le\tau_{k^s}(h).
\]
For $s=1$, equality holds.
If the assertion holds for $s-1$, then
\begin{align*}
\tau^s_k(h)
&=
\sum_{\substack{u_1,\ldots,u_k\ge1\\
u_1\cdots u_k=h}}
\tau^{s-1}_k(u_1\cdots u_k)\le
\sum_{\substack{u_1,\ldots,u_k\ge1\\
u_1\cdots u_k=h}}
\prod_{j=1}^{k}\tau^{s-1}_k(u_j)\le
\sum_{\substack{u_1,\ldots,u_k\ge1\\
u_1\cdots u_k=h}}
\prod_{j=1}^{k}\tau_{k^{s-1}}(u_j)
=
\tau_{k^s}(h).
\end{align*}

Summing over $h\le Z$ and first carrying out
the summation over the last factor, we obtain
\begin{align*}
\sum_{h\le Z}\tau_r(h)
&=
\sum_{\substack{u_1,\ldots,u_r\ge1\\
u_1\cdots u_r\le Z}}1=
\sum_{\substack{u_1,\ldots,u_{r-1}\ge1\\
u_1\cdots u_{r-1}\le Z}}
\left[
\frac{Z}{u_1\cdots u_{r-1}}
\right]\le
Z\left(\sum_{u\le Z}\frac{1}{u}\right)^{r-1}
\ll 
Z\bigl(\log(2Z)\bigr)^{r-1}.
\end{align*}
Consequently,
\begin{equation}\label{eq:divisor-third}
\sum_{h\le Z}\tau^s_k(h)
\le
\sum_{h\le Z}\tau_{k^s}(h)
\ll 
Z\bigl(\log(2Z)\bigr)^{k^s-1}.
\end{equation}

We return to the sum in~\eqref{eq:two-differences}.
Putting $Z=M^k$ in~\eqref{eq:divisor-third}
and using $\log(2M^k)\ll  l$, we have
\[
\sum_{h\le M^k}\tau^s_k(h)
\ll 
M^k l^{k^s-1}.
\]
Since $0\le v_h\le M$,
estimate~\eqref{eq:spacing} gives
\[
\sum_{h\le M^k}v_h^{\frac{s}{s-1}}
\le
M^{\frac{1}{s-1}}
\sum_{h\le M^k}v_h
\ll 
M^{n+\frac{1}{s-1}}Dl.
\]
Applying H\"older's inequality with exponents
$s$ and $\frac{s}{s-1}$, we obtain
\begin{align*}
\sum_{h\le M^k}\tau_k(h)v_h
&\le
\left(
\sum_{h\le M^k}\tau^s_k(h)
\right)^{\frac{1}{s}}
\left(
\sum_{h\le M^k}v_h^{\frac{s}{s-1}}
\right)^{\frac{s-1}{s}}\\
&\ll 
\left(M^k l^{k^s-1}\right)^{\frac{1}{s}}
\left(
M^{n+\frac{1}{s-1}}Dl
\right)^{\frac{s-1}{s}}=
M^n
D^{\frac{s-1}{s}}
l^{\frac{k^s+s-2}{s}}.
\end{align*}
Substituting this estimate
into~\eqref{eq:two-differences}
and using $k=n-1$, we obtain
\[
|T(\alpha,x,y)|^{2^k}
\ll 
M^{2^k-1}
+
M^{2^k}
D^{\frac{s-1}{s}}
l^{\frac{k^s+s-2}{s}}.
\]
Since $D\ge M^{-1}$ and $l\ge1$,
the first term is absorbed by the second.
Hence
\[
|T(\alpha,x,y)|
\ll 
MD^{\frac{s-1}{s\cdot 2^k}}
l^{\frac{k^s+s-2}{s\cdot2^k}}.
\]

Using $M\asymp y$,
$D\ll\frac{1}{q}+\frac{1}{y}+\frac{q}{y^n}$
and $l\ll\log(qy)$,
we obtain~\eqref{eq:log-weyl}.

We now prove~\eqref{eq:vmvt-log}.
Let $n\ge4$. For positive integers $s,k$, let $W_{s,k}(X)$
denote the number of integer solutions of the system
\[
\sum_{i=1}^{s}(u_i^j-v_i^j)=0,
\qquad j=1,\ldots,k,
\qquad 1\le u_i,v_i\le X.
\]

First suppose that $q\le M^n$.
By the estimate given
in~\cite[\S~8, proof of Theorem~1.5]{Wooley2012},
for every positive integer $s$ we have
\[
|T(\alpha;x,y)|
\ll
\left(
W_{s,n-1}(2M)
M^{\frac{n(n-1)}{2}}
\left(
\frac{1}{q}+\frac{1}{M}+\frac{q}{M^n}
\right)
\right)^{\frac{1}{2s}}
\log(2M).
\]
Put
\[
s_0=\frac{n(n-1)}{2}+1.
\]
By~\cite[Corollary~1.3]{Wooley2019},
with $k=n-1$, we have
\[
W_{s_0,n-1}(2M)
\ll
M^{2s_0-\frac{n(n-1)}{2}}.
\]
Substituting this estimate into the inequality for
$T(\alpha;x,y)$ with $s=s_0$, we obtain
\[
|T(\alpha;x,y)|
\ll
M\left(
\frac{1}{q}+\frac{1}{M}+\frac{q}{M^n}
\right)^{\frac{1}{n(n-1)+2}}
\log(2M).
\]
For $q>M^n$, the last inequality follows
from the trivial bound $|T(\alpha;x,y)|\le M$.

Finally, since $M\asymp y$ and
$\log(2M)\asymp\log y$ for $y\ge2$,
we obtain~\eqref{eq:vmvt-log}.

\section{Proof of Theorem \ref{thm:main}}
Let $N$ be sufficiently large, and let $H$
satisfy~\eqref{eq:Hrange}.

Put
\[
N_k=\mu_kN+H,\qquad k=1,2,3,
\qquad
H_3=N_3^{\frac{1}{n}}-(N_3-2H)^{\frac{1}{n}},
\]
and consider the sums
\[
V(\alpha,N_k,2H)
=
\sum_{N_k-2H<p\le N_k}e(\alpha p),
\qquad k=1,2,
\]
\[
T\left(\alpha,N_3^{\frac{1}{n}},H_3\right)
=
\sum_{(N_3-2H)^{\frac{1}{n}}<m\le N_3^{\frac{1}{n}}}
e(\alpha m^n).
\]
Taylor expansion gives
\begin{align}\label{H3}
N_3^{\frac{1}{n}}
&=
\mu_3^{\frac{1}{n}}N^{\frac{1}{n}}
\left(1+O\left(\frac{H}{N}\right)\right),
\nonumber\\
H_3
&=
\frac{2H}
{n\mu_3^{\frac{n-1}{n}}N^{\frac{n-1}{n}}}
\left(1+O\left(\frac{H^2}{N^2}\right)\right).
\end{align}
In particular,
\[
H_3\asymp\frac{H}{N^{\frac{n-1}{n}}},
\qquad
\frac{H_3}{N_3^{\frac{1}{n}}}
\ll\frac{H}{N}\le N^{-\delta}.
\]
Also put
\[
\tau
=
2n(n-1)N_3^{\frac{n-2}{n}}H_3
=
\frac{4(n-1)H}
{\mu_3^{\frac{1}{n}}N^{\frac{1}{n}}}
\left(1+O\left(\frac{H}{N}\right)\right),
\qquad
\kappa=\frac{1}{\tau}.
\]

By orthogonality, the number of solutions can be written
as
\begin{align}\label{eq:circle}
J(N,H)
=&
\int\limits_{-\kappa}^{1-\kappa}e(-\alpha N)
\bigl(V(\alpha,N_1,2H)+\theta_1(\alpha)\bigr)
\bigl(V(\alpha,N_2,2H)+\theta_2(\alpha)\bigr)
\nonumber\\
&\times
\left(
T\left(\alpha,N_3^{\frac{1}{n}},H_3\right)
+\theta_3(\alpha)
\right)\,d\alpha.
\end{align}
Here the corrections $\theta_k(\alpha)$ account for
possible terms at the lower endpoints of the intervals:
\[
\theta_k(\alpha)=
\begin{cases}
e\bigl(\alpha(N_k-2H)\bigr),
 &\text{if $N_k-2H\in\mathbb{P}$},\\
0,&\text{otherwise},
\end{cases}
\qquad k=1,2,
\]
and
\[
\theta_3(\alpha)=
\begin{cases}
e\bigl(\alpha(N_3-2H)\bigr),
 &\text{if $\sqrt[n]{N_3-2H}\in\mathbb N$},\\
0,&\text{otherwise}.
\end{cases}
\]
In particular, $|\theta_k(\alpha)|\le1$ for $k=1,2,3$.

Expanding the product in~\eqref{eq:circle}, we obtain
\begin{align}\label{eq:boundary-decomposition}
J(N,H)
=&
\int\limits_{-\kappa}^{1-\kappa}e(-\alpha N)
V(\alpha,N_1,2H)
V(\alpha,N_2,2H)
T\left(\alpha,N_3^{\frac{1}{n}},H_3\right)\,d\alpha
+R_1.
\end{align}
Since $|\theta_k(\alpha)|\le1$ for $k=1,2,3$,
the periodicity of the sums gives
\begin{align*}
|R_1|\le &
\int\limits_0^1
|V(\alpha,N_1,2H)V(\alpha,N_2,2H)|
\,d\alpha+\sum_{k=1}^2\int\limits_0^1
\left|
V(\alpha,N_k,2H)
T\left(\alpha,N_3^{\frac{1}{n}},H_3\right)
\right|\,d\alpha
\\
&+\sum_{k=1}^2\int\limits_0^1
|V(\alpha,N_k,2H)|\,d\alpha
+\int\limits_0^1
\left|T\left(\alpha,N_3^{\frac{1}{n}},H_3\right)\right|
\,d\alpha+1.
\end{align*}
By orthogonality and the
Brun--Titchmarsh inequality, we have
\[
\int\limits_0^1
|V(\alpha,N_k,2H)|^2\,d\alpha
=
\pi(N_k)-\pi(N_k-2H)
\ll\frac{H}{L},
\qquad k=1,2.
\]
Similarly,
\begin{align*}
\int\limits_0^1
\left|
T\left(\alpha,N_3^{\frac{1}{n}},H_3\right)
\right|^2\,d\alpha
&=
\left[N_3^{\frac{1}{n}}\right]
-
\left[(N_3-2H)^{\frac{1}{n}}\right]\le H_3+1
\ll \frac{H}{N^{\frac{n-1}{n}}}.
\end{align*}
By Cauchy's inequality,
\[
\int\limits_0^1|f(\alpha)g(\alpha)|\,d\alpha
\le
\left(
\int\limits_0^1|f(\alpha)|^2\,d\alpha
\right)^{\frac{1}{2}}
\left(
\int\limits_0^1|g(\alpha)|^2\,d\alpha
\right)^{\frac{1}{2}}.
\]
Applying this to the product of the two sums over primes
and using the estimates above, we obtain
\begin{align}
\int\limits_0^1
&\left|
V(\alpha,N_1,2H)
V(\alpha,N_2,2H)
\right|\,d\alpha
\le
\left(
\int\limits_0^1
|V(\alpha,N_1,2H)|^2\,d\alpha
\right)^{\frac{1}{2}}
\nonumber\\
&\times
\left(
\int\limits_0^1
|V(\alpha,N_2,2H)|^2\,d\alpha
\right)^{\frac{1}{2}}
\ll
\left(\frac{H}{L}\right)^{\frac{1}{2}}
\left(\frac{H}{L}\right)^{\frac{1}{2}}
=\frac{H}{L}.
\label{Est_int_primes}
\end{align}

Further, for $k=1,2$ we have
\begin{align*}
&\int\limits_0^1
\left|
V(\alpha,N_k,2H)
T\left(\alpha,N_3^{\frac{1}{n}},H_3\right)
\right|\,d\alpha
\le
\left(
\int\limits_0^1
|V(\alpha,N_k,2H)|^2\,d\alpha
\right)^{\frac{1}{2}}
\\
&\times
\left(
\int\limits_0^1
\left|T\left(\alpha,N_3^{\frac{1}{n}},H_3\right)\right|^2
\,d\alpha
\right)^{\frac{1}{2}}
\ll
\left(\frac{H}{L}\right)^{\frac{1}{2}}
\left(\frac{H}{N^{\frac{n-1}{n}}}\right)^{\frac{1}{2}}
=
\frac{H}{N^{\frac{n-1}{2n}}L^{\frac{1}{2}}}.
\end{align*}

For the integrals of the absolute value of a single sum, we apply
the same inequality, taking the second factor to be
one. For $k=1,2$, we obtain
\begin{align*}
\int\limits_0^1
|V(\alpha,N_k,2H)|\,d\alpha
&\le
\left(
\int\limits_0^1
|V(\alpha,N_k,2H)|^2\,d\alpha
\right)^{\frac{1}{2}}
\left(
\int\limits_0^1 1\,d\alpha
\right)^{\frac{1}{2}}\ll
\left(\frac{H}{L}\right)^{\frac{1}{2}}.
\end{align*}
Similarly,
\begin{align*}
&\int\limits_0^1
\left|T\left(\alpha,N_3^{\frac{1}{n}},H_3\right)\right|
\,d\alpha\le
\left(
\int\limits_0^1
\left|T\left(\alpha,N_3^{\frac{1}{n}},H_3\right)\right|^2
\,d\alpha
\right)^{\frac{1}{2}}
\left(
\int\limits_0^1 1\,d\alpha
\right)^{\frac{1}{2}}
\ll
\frac{H^{\frac{1}{2}}}{N^{\frac{n-1}{2n}}}.
\end{align*}

Substituting these estimates into the inequality for $R_1$,
we obtain
\[
|R_1|
\ll
\frac{H}{L}
+\frac{H}{N^{\frac{n-1}{2n}}L^{\frac{1}{2}}}
+\left(\frac{H}{L}\right)^{\frac{1}{2}}
+\frac{H^{\frac{1}{2}}}{N^{\frac{n-1}{2n}}}
+1.
\]
 
For sufficiently large $N$, we have
\[
\frac{H}{N^{\frac{n-1}{2n}}L^{\frac{1}{2}}}
\le \frac{H}{L},
\qquad
\frac{H^{\frac{1}{2}}}{N^{\frac{n-1}{2n}}}
\le
\left(\frac{H}{L}\right)^{\frac{1}{2}}
\le \frac{H}{L},
\qquad
1\le\frac{H}{L}.
\]
Consequently,
\[
|R_1|\ll\frac{H}{L}.
\]
Moreover, the condition on $H$ and the inequality
$\theta_n\le\frac{1}{2n}$ imply that
\[
\frac{N^{\frac{n-1}{n}}L^2}{H}
\le
N^{\theta_n-\frac{1}{n}}L^{2-A_n}
\le
N^{-\frac{1}{2n}}L^{2-A_n}
\rightarrow 0.
\]
Hence
\[
|R_1|
\ll
\frac{H}{L}
\ll
\frac{H^2}{N^{\frac{n-1}{n}}L^3}.
\]

Substituting the resulting estimate for $R_1$ into~\eqref{eq:boundary-decomposition}, we obtain


\begin{align}\label{eq:circle-reduced}
J(N,H)
=&
\int\limits_{-\kappa}^{1-\kappa}e(-\alpha N)
V(\alpha,N_1,2H)
V(\alpha,N_2,2H)
\nonumber\\
&\times
T\left(\alpha,N_3^{\frac{1}{n}},H_3\right)\,d\alpha
+
O\left(
\frac{H^2}{N^{\frac{n-1}{n}}L^3}
\right).
\end{align}

By Lemma~\ref{lemmaDirichle}, every number
$\alpha\in[-\kappa,1-\kappa)$ admits
a rational approximation
\[
\alpha=\frac{a}{q}+\lambda,
\qquad
(a,q)=1,
\qquad
1\le q\le\tau,
\qquad
|\lambda|\le\frac{1}{q\tau}.
\]
Let $Q=L^{B_n}$. For each reduced fraction
$\frac{a}{q}$ with $q\le Q$ and $0\le a\le q-1$,
define the interval
\[
\mathfrak M(q,a)=
\left\{
\alpha\in[-\kappa,1-\kappa):
\left|\alpha-\frac{a}{q}\right|
\le\frac{1}{q\tau}
\right\}.
\]

Write
\[
\mathfrak M
=
\bigcup_{q\le Q}
\ \bigcup_{\substack{0\le a\le q-1\\(a,q)=1}}
\mathfrak M(q,a),
\qquad
\mathfrak m=[-\kappa,1-\kappa)\setminus\mathfrak M.
\]

The intervals $\mathfrak M(q,a)$ are pairwise
disjoint.
Indeed, the distance between distinct
reduced fractions $\frac{a}{q}$
and $\frac{a_1}{q_1}$ is at least
$\frac{1}{qq_1}$, whereas
\[
\frac{1}{q\tau}+\frac{1}{q_1\tau}
=
\frac{q+q_1}{qq_1\tau}
<
\frac{1}{qq_1},
\]
since $q+q_1\le2Q<\tau$.

We divide $\mathfrak M$ into two sets:
\[
\mathfrak M_1
=
\left\{
\alpha\in\mathfrak M:
\left|\alpha-\frac{a}{q}\right|
\le\frac{L^2}{H}
\right\},
\]
\[
\mathfrak M_2
=
\left\{
\alpha\in\mathfrak M:
\frac{L^2}{H}
<
\left|\alpha-\frac{a}{q}\right|
\le\frac{1}{q\tau}
\right\}.
\]

Let $J(\mathfrak M_1)$,
$J(\mathfrak M_2)$, and $J(\mathfrak m)$
denote the integrals in~\eqref{eq:circle-reduced}
over the respective sets.
Then
\[
J(N,H)
=
J(\mathfrak M_1)+J(\mathfrak M_2)+J(\mathfrak m)
+
O\left(
\frac{H^2}{N^{\frac{n-1}{n}}L^3}
\right).
\]

{\bf Evaluation of the integral $J(\mathfrak M_1)$.}
Put
\[
F(\alpha)
=
V(\alpha,N_1,2H)
V(\alpha,N_2,2H)
T\left(\alpha,N_3^{\frac{1}{n}},H_3\right).
\]
We have
\[
J(\mathfrak M_1)
=
\sum_{q\le Q}
\sum_{\substack{0\le a\le q-1\\(a,q)=1}}
I(a,q),
\]
where
\begin{equation}\label{I(a,q)}
I(a,q)
=
e\left(-\frac{aN}{q}\right)
\int\limits_{-\frac{L^2}{H}}^{\frac{L^2}{H}}
F\left(\frac{a}{q}+\lambda\right)
e(-\lambda N)\,d\lambda.
\end{equation}
We apply Lemma~\ref{lem:primes-input},
putting
\[
x=N_k,\quad y=2H,\quad
A=[B_n]+2,\qquad b=[B_n]+1,
\quad \mathcal{L}=\log(N_kq),\quad k=1,2.
\]
For $q\le Q=L^{B_n}$, we have
\[
\mathcal L=L+O(\log L)\asymp L.
\]
All the hypotheses of this lemma hold
for sufficiently large $N$.
Indeed, since $[B_n]+1>B_n$
and $\log N_k=L+O(1)$, we have
\[
q\le Q=L^{B_n}
\le\mathcal{L}^{[B_n]+1}.
\]
Further, the condition on $H$ and the inequality
$1-\theta_n\ge\frac{5}{6}>\frac{5}{8}$ imply that
\[
2H\ge
N_k^{\frac{5}{8}}
\mathcal{L}^{\frac{7[B_n]+85}{4}}.
\]
Moreover,
\[
|\lambda|\le\frac{L^2}{H}
\le\frac{N_k}{8\pi H^2},
\]
since $N_k\asymp N$ and
\[
\frac{HL^2}{N}
\le N^{-\delta}L^2
\rightarrow0.
\]


Consequently,
\[
S_1\left(\frac{a}{q}+\lambda;N_k,2H\right)
=
\frac{\mu(q)}{\varphi(q)}
\frac{\sin(2\pi H\lambda)}{\pi\lambda}
e(\lambda\mu_kN)
+
O\left(\frac{H}{L^{[B_n]+2}}\right).
\]
We express $V(\alpha,N_k,2H)$ in terms of the sum
$S_1(\alpha,N_k,2H)$:
\begin{align*}
&V(\alpha,N_k,2H)
=
\frac{1}{\log(\mu_kN)}
\sum_{N_k-2H<p\le N_k}
\log(\mu_kN)e(\alpha p)
\\
&=
\frac{1}{\log(\mu_kN)}
\sum_{N_k-2H<p\le N_k}
(\log p)e(\alpha p)+
\frac{1}{\log(\mu_kN)}
\sum_{N_k-2H<p\le N_k}
\log\left(\frac{\mu_kN}{p}\right)e(\alpha p)
\\
&=
\frac{S_1(\alpha,N_k,2H)}{\log(\mu_kN)}
+
\frac{1}{\log(\mu_kN)}
\sum_{N_k-2H<p\le N_k}
\log\left(\frac{\mu_kN}{p}\right)e(\alpha p)
\\
&\quad-
\frac{1}{\log(\mu_kN)}
\sum_{r\ge2}
\ \sum_{N_k-2H<p^r\le N_k}
(\log p)e(\alpha p^r).
\end{align*}

We estimate the last two terms.
In the summation interval,
$|p-\mu_kN|\le H$, so
\[
\log\left(\frac{\mu_kN}{p}\right)
=
-\log\left(
1+\frac{p-\mu_kN}{\mu_kN}
\right)
\ll\frac{H}{N}.
\]
Hence the Brun--Titchmarsh inequality gives
\begin{align*}
&\frac{1}{\log(\mu_kN)}
\left|
\sum_{N_k-2H<p\le N_k}
\log\left(\frac{\mu_kN}{p}\right)e(\alpha p)
\right|\ll
\frac{H}{NL}
\sum_{N_k-2H<p\le N_k}1
\ll\frac{H^2}{NL^2}.
\end{align*}

For the term involving prime powers, we have
\begin{align*}
&\frac{1}{\log(\mu_kN)}
\left|
\sum_{r\ge2}
\ \sum_{N_k-2H<p^r\le N_k}
(\log p)e(\alpha p^r)
\right|\ll
\frac{1}{L}
\sum_{2\le r\le\frac{\log N_k}{\log2}}
\ \sum_{p\le N_k^{\frac{1}{r}}}\log p
\\
&\qquad\le
\frac{\log N_k}{L}
\sum_{2\le r\le\frac{\log N_k}{\log2}}
N_k^{\frac{1}{r}}
\ll N^{\frac{1}{2}}L.
\end{align*}

Thus,
\[
V(\alpha,N_k,2H)
=
\frac{S_1(\alpha,N_k,2H)}{\log(\mu_kN)}
+
O\left(
\frac{H^2}{NL^2}+N^{\frac{1}{2}}L
\right).
\]


We show that the additional error terms are absorbed
by $\frac{H}{L^{A+1}}$, where $A=[B_n]+2$.
Indeed, the upper bound for $H$ implies that
\[
\frac{H^2}{NL^2}\cdot\frac{L^{A+1}}{H}
=
\frac{HL^{A-1}}{N}
\le
N^{-\delta}L^{A-1}
\rightarrow0.
\]
The lower bound for $H$ gives
\[
N^{\frac{1}{2}}L\cdot\frac{L^{A+1}}{H}
=
\frac{N^{\frac{1}{2}}L^{A+2}}{H}
\le
N^{\theta_n-\frac{1}{2}}L^{A+2-A_n}
\rightarrow0,
\]
since $\theta_n\le\frac{1}{6}$.

Substituting the formula for $S_1$ obtained above
and using $\log(\mu_kN)\asymp L$, we obtain
\begin{align*}
V\left(\frac{a}{q}+\lambda,N_k,2H\right)
={}&
\frac{\mu(q)}{\varphi(q)}
\frac{\sin(2\pi H\lambda)}{\pi\lambda}
\frac{e(\lambda\mu_kN)}{\log(\mu_kN)}+
O\left(\frac{H}{L^{A+1}}\right),
\qquad k=1,2.
\end{align*}
Multiplying the formulas for $k=1,2$, we obtain
\begin{align*}
&V\left(\frac{a}{q}+\lambda,N_1,2H\right)
V\left(\frac{a}{q}+\lambda,N_2,2H\right)
\\
&\quad=
\frac{\mu^2(q)}{\varphi^2(q)}
\frac{\sin^2(2\pi H\lambda)}{\pi^2\lambda^2}
\frac{e\bigl(\lambda(\mu_1+\mu_2)N\bigr)}
{\log(\mu_1N)\log(\mu_2N)}
+
O\left(
\frac{H^2}{\varphi(q)L^{A+2}}
\right).
\end{align*}
The product of the error terms is absorbed by the stated error term,
since
$$
\varphi(q)L^{-A}\le L^{B_n-A}\ll 1.
$$
For the sum $T$, we use Lemma~\ref{exposumm1},
putting
\[
x=N_3^{\frac{1}{n}},
\qquad y=H_3.
\]
For sufficiently large $N$, we have
\[
\frac{H_3}{N_3^{\frac{1}{n}}}
\ll\frac{H}{N}
\le\frac{1}{N^\delta}
\rightarrow 0,
\]
\[
q\le Q<\tau
=
2n(n-1)N_3^{\frac{n-2}{n}}H_3.
\]
Hence the conditions $1\le y\le 0.01x$,
$q\le\tau$, and $\tau\ge2n(n-1)x^{n-2}y$
are satisfied.

We check the restriction on $\lambda$.
On the set $\mathfrak M_1$, we obtain
\begin{align*}
2nqN_3^{\frac{n-1}{n}}|\lambda|
&\le
\frac{2nqN_3^{\frac{n-1}{n}}L^2}{H}
\le
\frac{
2nN_3^{\frac{n-1}{n}}L^{B_n+2}
}{
N^{1-\theta_n}L^{A_n}
}
\ll
\frac{L^{B_n+2-A_n}}{N^{\frac{1}{2n}}}
\rightarrow 0.
\end{align*}
Consequently, for sufficiently large $N$,
\[
|\lambda|
\le
\frac{1}{2nqN_3^{\frac{n-1}{n}}}.
\]

Now Lemma~\ref{exposumm1} gives
\begin{equation*}
T\left(
\frac{a}{q}+\lambda,N_3^{\frac{1}{n}},H_3
\right)
=\frac{H_3}{q}S(a,q)
\gamma\left(\lambda;N_3^{\frac{1}{n}},H_3\right)+
O\left(q^{\frac{1}{2}+\varepsilon}\right),
\end{equation*}
where
\[
\gamma\left(\lambda;N_3^{\frac{1}{n}},H_3\right)
=
\int\limits_{-\frac{1}{2}}^{\frac{1}{2}}
e\left(
\lambda
\left(
N_3^{\frac{1}{n}}
+
H_3\left(t-\frac{1}{2}\right)
\right)^n
\right)\,dt.
\]

By the definition of $H_3$,
\[
2H
=
N_3-
\left(N_3^{\frac{1}{n}}-H_3\right)^n
=
nN_3^{\frac{n-1}{n}}H_3
+
O\left(N_3^{\frac{n-2}{n}}H_3^2\right).
\]
Consequently, for
$-\frac{1}{2}\le t\le\frac{1}{2}$
\begin{align*}
&\left(
N_3^{\frac{1}{n}}+
H_3\left(t-\frac{1}{2}\right)
\right)^n=
N_3+
nN_3^{\frac{n-1}{n}}H_3
\left(t-\frac{1}{2}\right)
+
O\left(N_3^{\frac{n-2}{n}}H_3^2\right)\\
&\qquad=
\mu_3N+2Ht+R_2(t),
\end{align*}
where
\[
|R_2(t)|
\ll
N_3^{\frac{n-2}{n}}H_3^2
\ll N^{\frac{n-2}{n}}\left(\frac{H}{N^{\frac{n-1}{n}}}\right)^2
\ll\frac{H^2}{N}.
\]
Using the inequality $|e(u)-1|\le2\pi|u|$,
we obtain
\begin{align*}
&\gamma\left(\lambda;N_3^{\frac{1}{n}},H_3\right)
=e(\lambda\mu_3N)
\int\limits_{-\frac{1}{2}}^{\frac{1}{2}}
e(2\lambda Ht)e(\lambda R_2(t))\,dt=e(\lambda\mu_3N)
\int\limits_{-\frac{1}{2}}^{\frac{1}{2}}
e(2\lambda Ht)dt\\
&+e(\lambda\mu_3N)
\int\limits_{-\frac{1}{2}}^{\frac{1}{2}}
e(2\lambda Ht)\Big(e(\lambda R_2(t))-1\Big)\,dt=
e(\lambda\mu_3N)
\frac{\sin(2\pi H\lambda)}{2\pi H\lambda}
+
O\left(\frac{|\lambda|H^2}{N}\right).
\end{align*}
Thus,
\begin{align*}
&T\left(
\frac{a}{q}+\lambda,N_3^{\frac{1}{n}},H_3
\right)=
\frac{H_3S(a,q)}{q}
\frac{\sin(2\pi H\lambda)}{2\pi H\lambda}
e(\lambda\mu_3N)
+R_3,
\end{align*}
where
\[
R_3
=
O\left(
\frac{H_3|S(a,q)|}{q}
\frac{|\lambda|H^2}{N}
\right)
+
O\left(q^{\frac{1}{2}+\varepsilon}\right).
\]
By Lemma~\ref{orsumm} and the condition
$|\lambda|\le\frac{L^2}{H}$, we have
\begin{equation}\label{eq:R3-bound}
|R_3|
\ll
\frac{H_3HL^2}{q^{\frac{1}{n}}N}
+q^{\frac{1}{2}+\varepsilon}\ll
\frac{H^2L^2}
{q^{\frac{1}{n}}N^{\frac{2n-1}{n}}}
+q^{\frac{1}{2}+\varepsilon}.
\end{equation}
Using $\mu_1+\mu_2+\mu_3=1$, we obtain
\[
F(\alpha)
=
\frac{H_3}{2H\log(\mu_1N)\log(\mu_2N)}
\frac{\mu^2(q)S(a,q)}{q\varphi^2(q)}
\frac{\sin^3(2\pi H\lambda)}{\pi^3\lambda^3}
e(\lambda N)
+
R_4,
\]
where
\begin{align*}
R_4={}&
\frac{\mu^2(q)}{\varphi^2(q)}
\frac{\sin^2(2\pi H\lambda)}{\pi^2\lambda^2}
\frac{e\bigl(\lambda(\mu_1+\mu_2)N\bigr)}
{\log(\mu_1N)\log(\mu_2N)}
R_3
\\
&+
\frac{H_3S(a,q)}{q}
\frac{\sin(2\pi H\lambda)}{2\pi H\lambda}
e(\lambda\mu_3N)
O\left(\frac{H^2}{\varphi(q)L^{A+2}}\right)+
R_3O\left(\frac{H^2}{\varphi(q)L^{A+2}}\right).
\end{align*}

By Lemma~\ref{orsumm} and the inequality
$|\sin u|\le|u|$, we obtain
\begin{align*}
|R_4|
\ll{}&
\frac{H^2}{\varphi^2(q)L^2}|R_3|
+
\frac{H^2H_3}
{q^{\frac{1}{n}}\varphi(q)L^{A+2}}+
\frac{H^2}{\varphi(q)L^{A+2}}|R_3|.
\end{align*}
The last term is absorbed by the first,
since $\varphi(q)L^{-A}\le L^{B_n-A}\ll1$.
Hence
\[
|R_4|\ll
\frac{H^2}{\varphi^2(q)L^2}|R_3|
+\frac{H^2H_3}
{q^{\frac{1}{n}}\varphi(q)L^{A+2}}.
\]

Substituting estimate~\eqref{eq:R3-bound} for $R_3$ and using
$H_3\ll\frac{H}{N^{\frac{n-1}{n}}}$, we obtain
\begin{align*}
|R_4|
\ll{}&
\frac{H^3}
{q^{\frac{1}{n}}\varphi(q)
N^{\frac{n-1}{n}}L^{A+2}}+
\frac{H^4}
{q^{\frac{1}{n}}\varphi^2(q)
N^{\frac{2n-1}{n}}}+
\frac{H^2q^{\frac{1}{2}+\varepsilon}}
{\varphi^2(q)L^2}.
\end{align*}
The ratio of the second term to the first is
\[
\frac{HL^{A+2}}{\varphi(q)N}
\le
N^{-\delta}L^{A+2}
\rightarrow0.
\]
Consequently,
\begin{equation}\label{eq:R4-bound}
|R_4|
\ll
\frac{H^3}
{q^{\frac{1}{n}}\varphi(q)
N^{\frac{n-1}{n}}L^{A+2}}
+
\frac{H^2q^{\frac{1}{2}+\varepsilon}}
{\varphi^2(q)L^2}.
\end{equation}

Substituting the expression for $F$
into~\eqref{I(a,q)}, we have
\begin{equation}\label{eq:Iaq-KH}
I(a,q)
=
\frac{H_3}
{2H\log(\mu_1N)\log(\mu_2N)}
\frac{\mu^2(q)S(a,q)}{q\varphi^2(q)}
e\left(-\frac{aN}{q}\right)K(H)
+
R_5(a,q),
\end{equation}
where
\[
K(H)
=
\int\limits_{-\frac{L^2}{H}}^{\frac{L^2}{H}}
\frac{\sin^3(2\pi H\lambda)}
{\pi^3\lambda^3}\,d\lambda,
\]
\[
R_5(a,q)
=
e\left(-\frac{aN}{q}\right)
\int\limits_{-\frac{L^2}{H}}^{\frac{L^2}{H}}
R_4e(-\lambda N)\,d\lambda.
\]
It follows from~\eqref{eq:R4-bound} that
\begin{align}\label{eq:R5-bound}
|R_5(a,q)|
&\le
\int\limits_{-\frac{L^2}{H}}^{\frac{L^2}{H}}
|R_4|\,d\lambda\ll
\frac{L^2}{H}
\left(
\frac{H^3}
{q^{\frac{1}{n}}\varphi(q)
N^{\frac{n-1}{n}}L^{A+2}}
+
\frac{H^2q^{\frac{1}{2}+\varepsilon}}
{\varphi^2(q)L^2}
\right)
\nonumber\\
&=
\frac{H^2}
{q^{\frac{1}{n}}\varphi(q)
N^{\frac{n-1}{n}}L^A}
+
\frac{Hq^{\frac{1}{2}+\varepsilon}}
{\varphi^2(q)}.
\end{align}

Making the substitution $u=2\pi H\lambda$
and applying Lemma~\ref{vatson}, we obtain
\[
K(H)
=
\frac{8H^2}{\pi}
\int\limits_0^{2\pi L^2}\frac{\sin^3u}{u^3}\,du,
\qquad
\int\limits_0^\infty\frac{\sin^3u}{u^3}\,du
=
\frac{3\pi}{8}.
\]
Estimating the tail of the integral,
\[
\left|
\int\limits_{2\pi L^2}^{\infty}
\frac{\sin^3u}{u^3}\,du
\right|
\le
\int\limits_{2\pi L^2}^{\infty}\frac{du}{u^3}
\ll L^{-4},
\]
we obtain
\begin{equation}\label{eq:KH-asymptotic}
K(H)=3H^2+O\left(\frac{H^2}{L^4}\right).
\end{equation}
Thus, substituting the value of $K(H)$ from \eqref{eq:KH-asymptotic} into \eqref{eq:Iaq-KH}, we obtain
\begin{equation}\label{eq:Iaq-asymptotic}
I(a,q)
=
\frac{3HH_3}
{2\log(\mu_1N)\log(\mu_2N)}
\frac{\mu^2(q)S(a,q)}{q\varphi^2(q)}
e\left(-\frac{aN}{q}\right)
+
R_6(a,q),
\end{equation}
where~\eqref{eq:R5-bound} gives
\begin{align}\label{eq:R6-bound}
|R_6(a,q)|
&\ll
|R_5(a,q)|
+
\frac{HH_3}
{q^{\frac{1}{n}}\varphi^2(q)L^6}
\nonumber\\
&\ll
\frac{H^2}
{q^{\frac{1}{n}}\varphi(q)
N^{\frac{n-1}{n}}L^A}
+
\frac{Hq^{\frac{1}{2}+\varepsilon}}
{\varphi^2(q)}
+
\frac{H^2}
{q^{\frac{1}{n}}\varphi^2(q)
N^{\frac{n-1}{n}}L^6}.
\end{align}

Summing~\eqref{eq:Iaq-asymptotic}
over the central arcs, we obtain
\begin{equation}\label{eq:JM1-R7}
J(\mathfrak M_1)
=\frac{3HH_3}
{2\log(\mu_1N)\log(\mu_2N)}\sum_{q\le Q}
\frac{\mu^2(q)}{q\varphi^2(q)}
\sum_{\substack{a=0\\(a,q)=1}}^{q-1}
S(a,q)e\left(-\frac{aN}{q}\right)
+
R_7,
\end{equation}
where
\[
R_7
=
\sum_{q\le Q}
\sum_{\substack{0\le a\le q-1\\(a,q)=1}}
R_6(a,q).
\]

For each $q$, the number of terms in the inner
sum is $\varphi(q)$.
Hence~\eqref{eq:R6-bound} implies that
\begin{equation*}
|R_7|
\ll
\frac{H^2}{N^{\frac{n-1}{n}}L^A}
\sum_{q\le Q}\frac{1}{q^{\frac{1}{n}}}
+H\sum_{q\le Q}
\frac{q^{\frac{1}{2}+\varepsilon}}{\varphi(q)}
+\frac{H^2}{N^{\frac{n-1}{n}}L^6}
\sum_{q\le Q}
\frac{1}{q^{\frac{1}{n}}\varphi(q)}.
\end{equation*}

We use the identity
\[
\frac{q}{\varphi(q)}
=
\sum_{d\mid q}\frac{\mu^2(d)}{\varphi(d)}.
\]
Changing the order of summation and putting $q=dm$,
we obtain
\begin{align*}
\sum_{q\le Q}
\frac{1}{q^{\frac{1}{n}}\varphi(q)}
=\sum_{q\le Q}
\frac{1}{q^{1+\frac{1}{n}}}
\sum_{d\mid q}\frac{\mu^2(d)}{\varphi(d)}=
&\sum_{d\le Q}
\frac{\mu^2(d)}
{d^{1+\frac{1}{n}}\varphi(d)}
\sum_{m\le Q/d}\frac{1}{m^{1+\frac{1}{n}}}\le
\left(
\sum_{m=1}^{\infty}
\frac{1}{m^{1+\frac{1}{n}}}
\right)^2
\ll 1.
\end{align*}
For the remaining sums, we have
\[
\begin{aligned}
\sum_{q\le Q}\frac{1}{q^{\frac{1}{n}}}
&\le
1+\int\limits_1^Q t^{-\frac{1}{n}}\,dt=
1+\frac{n}{n-1}
\left(Q^{\frac{n-1}{n}}-1\right)\le
\frac{n}{n-1}Q^{\frac{n-1}{n}},
\end{aligned}
\]
\[
\sum_{q\le Q}
\frac{q^{\frac{1}{2}+\varepsilon}}{\varphi(q)}
\le
\sum_{q\le Q}q^{\frac{1}{2}+\varepsilon}
\ll Q^{\frac{3}{2}+\varepsilon}.
\]
Using $Q=L^{B_n}$, we obtain
\[
|R_7|
\ll
\frac{H^2}
{N^{\frac{n-1}{n}}
L^{A-\frac{n-1}{n}B_n}}
+
HQ^{\frac{3}{2}+\varepsilon}
+
\frac{H^2}{N^{\frac{n-1}{n}}L^6}.
\]

The definition of $B_n$ gives $B_n>2n$.
Since $A=[B_n]+2>B_n+1$, we have
\[
A-\frac{n-1}{n}B_n
>
1+\frac{B_n}{n}
>
3.
\]
For the second term, the condition on $H$ gives
\begin{align*}
HQ^{\frac{3}{2}+\varepsilon}=\frac{H^2}{N^{\frac{n-1}{n}}L^3}\frac{N^{\frac{n-1}{n}}L^3
Q^{\frac{3}{2}+\varepsilon}}{H}
&\le
\frac{H^2}{N^{\frac{n-1}{n}}L^3}N^{-\frac{1}{2n}}
L^{(\frac{3}{2}+\varepsilon)B_n+3-A_n}\ll \frac{H^2}{N^{\frac{n-1}{n}}L^3}.
\end{align*}
Thus,
\begin{equation}\label{eq:R7-bound}
|R_7|
\ll
\frac{H^2}{N^{\frac{n-1}{n}}L^3}.
\end{equation}

It remains to replace the sum over $q$
by an infinite series.
Write
\[
\Phi(q,N)
=
\sum_{\substack{0\le a\le q-1\\(a,q)=1}}
S(a,q)e\left(-\frac{aN}{q}\right).
\]
By the Chinese remainder theorem, the function
$\Phi(q,N)$ is multiplicative in $q$:
\[
\Phi(q_1q_2,N)=\Phi(q_1,N)\Phi(q_2,N),
\qquad (q_1,q_2)=1.
\]

For a prime $p$, we have
\begin{align*}
\Phi(p,N)
&=
\sum_{u=1}^{p}\sum_{a=1}^{p-1}
e\left(\frac{a(u^n-N)}{p}\right)
=(p-1)\rho(N,p)
-\bigl(p-\rho(N,p)\bigr)
\\
&=
p\bigl(\rho(N,p)-1\bigr),
\end{align*}
where $\rho(N,p)$ is the number of solutions of the congruence
\[
u^n\equiv N\pmod p.
\]
Since $0\le\rho(N,p)\le n$,
multiplicativity gives
\begin{align*}
\sum_{q=1}^{\infty}
\frac{\mu^2(q)|\Phi(q,N)|}
{q^{\frac{1}{2}}\varphi^2(q)}
&=
\prod_p
\left(
1+
\frac{p^{\frac{1}{2}}|\rho(N,p)-1|}
{(p-1)^2}
\right)\le
\prod_p
\left(
1+\frac{(n-1)p^{\frac{1}{2}}}{(p-1)^2}
\right) \\ 
& \le \prod_p
\left(
1+\frac{1}{p^{\frac{3}{2}}}
\right)^{4(n-1)}=\left(\sum_{m=1}^{\infty}\frac{\mu^2(m)}{m^{\frac{3}{2}}}\right)^{4(n-1)}\ll 1.
\end{align*}
Consequently, the series
\[
\mathfrak S(N)
=
\sum_{q=1}^{\infty}
\frac{\mu^2(q)\Phi(q,N)}{q\varphi^2(q)}
\]
converges absolutely and, by multiplicativity,
can be written as
\[
\mathfrak S(N)
=
\prod_p
\left(
1+\frac{\rho(N,p)-1}{(p-1)^2}
\right).
\]
For $p=2$, we have $\rho(N,2)=1$.
For $p\ge3$, the inequalities $0\le\rho(N,p)\le n$ imply that
\[
1-\frac{1}{(p-1)^2}
\le
1+\frac{\rho(N,p)-1}{(p-1)^2}
\le
1+\frac{n-1}{(p-1)^2}.
\]
Since the series
$\sum\limits_{p\ge3}\frac{1}{(p-1)^2}$ converges,
we obtain
\[
0<c=
\prod_{p\ge3}
\left(1-\frac{1}{(p-1)^2}\right)
\le\mathfrak S(N),
\]
where $c$ is an absolute constant, and
\[
\mathfrak S(N)
\le
\prod_{p\ge3}
\left(1+\frac{n-1}{(p-1)^2}\right)
=C_n<\infty.
\]

We estimate the tail of the series. Using the bound obtained above
for the weighted absolute series, we have
\begin{align}\label{eq:singular-series-tail}
\left|
\sum_{q>Q}
\frac{\mu^2(q)\Phi(q,N)}{q\varphi^2(q)}
\right|
&\le
Q^{-\frac{1}{2}}
\sum_{q>Q}
\frac{\mu^2(q)|\Phi(q,N)|}
{q^{\frac{1}{2}}\varphi^2(q)}\ll
Q^{-\frac{1}{2}}
=
L^{-\frac{B_n}{2}}
\ll L^{-2},
\end{align}
since $B_n>4$.

Thus, in view of \eqref{eq:singular-series-tail}, we have
\[
\sum_{q\le Q}
\frac{\mu^2(q)\Phi(q,N)}{q\varphi^2(q)}
=
\mathfrak S(N)
+
O\left(\frac{1}{L^2}\right).
\]

Substituting this identity
into~\eqref{eq:JM1-R7}, we see that the error in replacing the finite
sum by the infinite series is bounded by
\[
\frac{HH_3}{L^2}\cdot\frac{1}{L^2}
=
\frac{HH_3}{L^4}
\ll
\frac{H^2}{N^{\frac{n-1}{n}}L^4}.
\]
Hence, using~\eqref{eq:R7-bound}, we obtain
\[
J(\mathfrak M_1)
=
\frac{3HH_3\mathfrak S(N)}
{2\log(\mu_1N)\log(\mu_2N)}
+
O\left(
\frac{H^2}{N^{\frac{n-1}{n}}L^3}
\right).
\]

Using~\eqref{H3}, the relation
\[
\frac{1}{\log(\mu_1N)\log(\mu_2N)}
=
\frac{1}{L^2}+O\left(\frac{1}{L^3}\right)
\]
and the inequalities $\frac{H}{N}\le N^{-\delta}$,
$\mathfrak S(N)\ll1$, we obtain

\[
J(\mathfrak M_1)
=
\frac{3\mathfrak S(N)H^2}
{n\mu_3^{\frac{n-1}{n}}
N^{\frac{n-1}{n}}L^2}
+
O\left(
\frac{H^2}{N^{\frac{n-1}{n}}L^3}
\right).
\]

{\bf Estimation of $J(\mathfrak M_2)$.}
In the subsequent estimates, put
\[
x=N_3^{\frac{1}{n}},
\qquad
y=H_3.
\]
For $\alpha\in\mathfrak M_2$, we have
\[
\alpha=\frac{a}{q}+\lambda,
\qquad
q\le Q,\quad(a,q)=1,
\qquad
\frac{L^2}{H}<|\lambda|\le\frac{1}{q\tau}.
\]

First suppose that
\[
\frac{L^2}{H}<|\lambda|
\le\frac{1}{2nqx^{n-1}}.
\]
By Lemma~\ref{exposumm1},
\[
T(\alpha,x,y)
=
\frac{S(a,q)}{q}
\int\limits_{x-y}^{x}e(\lambda u^n)\,du
+
O\left(q^{\frac{1}{2}+\varepsilon}\right).
\]
Integrating by parts and using $y=o(x)$,
we obtain
\[
\int\limits_{x-y}^{x}
e(\lambda u^n)\,du
=\left.
\frac{e(\lambda u^n)}{2n\pi i\lambda u^{n-1}}
\right|_{x-y}^{x}+
\frac{n-1}{2n\pi i\lambda}
\int\limits_{x-y}^{x}
\frac{e(\lambda u^n)}{u^n}\,du.
\]
Consequently, since $x-y\asymp x$, we have
\[
\left|\ 
\int\limits_{x-y}^{x}
e(\lambda u^n)\,du
\right|
\ll
\frac{1}{|\lambda|}
\left(
\frac{1}{(x-y)^{n-1}}
+\frac{1}{x^{n-1}}
+\int\limits_{x-y}^{x}
\frac{du}{u^n}
\right)\ll
\frac{1}{|\lambda|x^{n-1}}.
\]
Hence, by Lemma~\ref{orsumm} and the relation $x^{n-1}y\asymp H$, we have
\[
|T(\alpha,x,y)|
\ll
\frac{1}{q^{\frac{1}{n}}|\lambda|x^{n-1}}
+
q^{\frac{1}{2}+\varepsilon}
\ll
\frac{y}{L^2}+Q^{\frac{1}{2}+\varepsilon}
\ll\frac{y}{L^2}.
\]
Now suppose that
\[
\frac{1}{2nqx^{n-1}}
<|\lambda|\le\frac{1}{q\tau}.
\]
By Lemma~\ref{exposumm2},
\[
|T(\alpha,x,y)|
\ll
q^{1-\frac{1}{n}}\log q
+
\min\left(
yq^{-\frac{1}{n}},
x^{\frac{1}{2}}q^{\frac{1}{2}-\frac{1}{n}}
\right).
\]
Using the inequality
\[
\min(A,B)\le A^{1-\frac{2}{n}}B^{\frac{2}{n}},
\qquad A,B>0,
\]
we obtain
\begin{align*}
&\min\left(
yq^{-\frac{1}{n}},
x^{\frac{1}{2}}q^{\frac{1}{2}-\frac{1}{n}}
\right)\le
\left(yq^{-\frac{1}{n}}\right)^{1-\frac{2}{n}}
\left(
x^{\frac{1}{2}}q^{\frac{1}{2}-\frac{1}{n}}
\right)^{\frac{2}{n}}
=
x^{\frac{1}{n}}y^{1-\frac{2}{n}}.
\end{align*}
Hence
\[
|T(\alpha,x,y)|
\ll
Q^{1-\frac{1}{n}}\log Q
+
x^{\frac{1}{n}}y^{1-\frac{2}{n}}.
\]
The first term is of smaller order
than $yL^{-2}$.

For the second term, we have
\begin{align*}
x^{\frac{1}{n}}y^{1-\frac{2}{n}}
&=
\frac{y}{L^2}
\left(\frac{x^{\frac{1}{2}}L^n}{y}\right)^{\frac{2}{n}}\ll
\frac{y}{L^2}
\left(
\frac{N^{1-\frac{1}{2n}}L^n}{H}
\right)^{\frac{2}{n}}\le
\frac{y}{L^2}
\left(
N^{\theta_n-\frac{1}{2n}}L^{n-A_n}
\right)^{\frac{2}{n}}
\ll
\frac{y}{L^2},
\end{align*}
since $\theta_n\le\frac{1}{2n}$ and $A_n\ge n$.

Thus,
\[
\sup_{\alpha\in\mathfrak M_2}|T(\alpha,x,y)|
\ll\frac{y}{L^2}.
\]
Applying~\eqref{Est_int_primes}, we obtain
\begin{align}\label{eq:JM2-bound}
|J(\mathfrak M_2)|
&\le
\sup_{\alpha\in\mathfrak M_2}|T(\alpha,x,y)|
\int\limits_{\mathfrak M_2}
\left|
V(\alpha,N_1,2H)
V(\alpha,N_2,2H)
\right|\,d\alpha
\nonumber\\
&\ll
\frac{y}{L^2}\cdot\frac{H}{L}
\ll
\frac{H^2}{N^{\frac{n-1}{n}}L^3}.
\end{align}

{\bf Estimation of $J(\mathfrak m)$.}
Put
\[
(\sigma_n,\gamma_n)=
\begin{cases}
\left(\frac{1}{6},\frac{3}{4}\right), & n=3,\\[4pt]
\left(\frac{1}{12},\frac{7}{6}\right), & n=4,\\[4pt]
\left(\frac{1}{n(n-1)+2},1\right), & n\ge5.
\end{cases}
\]
Then $B_n\sigma_n=2+\gamma_n$.

We apply Lemma~\ref{lemmaDirichle} again,
putting
\[
P=\frac{y^n}{Q}.
\]
For each $\alpha\in\mathfrak m$,
there exist integers $b,r$ such that
\[
\alpha=\frac{b}{r}+\beta,
\qquad
(b,r)=1,
\qquad
1\le r\le P,
\qquad
|\beta|\le\frac{1}{rP}\le\frac{1}{r^2}.
\]
We consider the cases {\bf 1.} $r>Q$ and {\bf 2.} $r\le Q$.

{\bf Case 1.}
Using estimate~\eqref{eq:log-weyl}
for $n=3$ and $4$, $s=3$,
and estimate~\eqref{eq:vmvt-log} for $n\ge5$,
we obtain, for $Q<r\le P$,
\begin{align*}
|T(\alpha,x,y)|
&\ll
y\left(
\frac{1}{r}+\frac{1}{y}+\frac{r}{y^n}
\right)^{\sigma_n}
\bigl(\log(ry)\bigr)^{\gamma_n}\\
&\ll
yQ^{-\sigma_n}L^{\gamma_n}
=
yL^{-B_n\sigma_n+\gamma_n}
=
\frac{y}{L^2},
\end{align*}
since $Q<y, r\le P=\frac{y^n}{Q}$ and $\log(ry)\ll L$.

{\bf Case 2.}
Let $r\le Q$.
Since $\alpha\in\mathfrak m$, we have
$|\beta|>\frac{1}{r\tau}.$

Consequently,
\[
\frac{1}{2n(n-1)r x^{n-2}y}
<|\beta|
\le\frac{Q}{ry^n}
\le\frac{1}{ry^{n-1}}.
\]
We apply the estimate
of Lemma~\ref{ShortExposumIR}.
Using both bounds for $|\beta|$, we obtain
\begin{align*}
|T(\alpha,x,y)|
\ll&
r^{1-\frac{1}{n}}L
+r^{\frac{1}{2}-\frac{1}{n}}y^{\frac{1}{2}}+
r^{\frac{1}{2}-\frac{1}{n}}Q^{\frac{1}{2}}
x^{\frac{n-2}{2}}y^{1-\frac{n}{2}}\\
\ll{}&
Q^{1-\frac{1}{n}}L
+Q^{\frac{1}{2}-\frac{1}{n}}y^{\frac{1}{2}}+
Q^{1-\frac{1}{n}}
x^{\frac{n-2}{2}}y^{1-\frac{n}{2}}.
\end{align*}

Since
\[
y\gg N^{\frac{1}{n}-\theta_n}L^{A_n},
\qquad x\asymp N^{\frac{1}{n}},
\]
we have
\[
Q^{1-\frac{1}{n}}L=\frac{y}{L^2}\cdot \frac{Q^{1-\frac{1}{n}}L^3}{y}
\ll\frac{y}{L^2}\cdot
N^{\theta_n-\frac{1}{n}}
L^{\frac{n-1}{n}B_n+3-A_n}
\ll \frac{y}{L^2}
\]
and
\[
Q^{\frac{1}{2}-\frac{1}{n}}y^{\frac{1}{2}}=\frac{y}{L^2}\cdot\frac{Q^{\frac{1}{2}-\frac{1}{n}}L^2}{y^{\frac{1}{2}}}
\ll\frac{y}{L^2}\cdot
N^{\frac{\theta_n}{2}-\frac{1}{2n}}
L^{\frac{n-2}{2n}B_n+2-\frac{A_n}{2}}
\ll\frac{y}{L^2}.
\]
For the last term, we obtain
\begin{align*}
Q^{1-\frac{1}{n}}
x^{\frac{n-2}{2}}y^{1-\frac{n}{2}}&=\frac{y}{L^2}\cdot\frac{Q^{1-\frac{1}{n}}
x^{\frac{n-2}{2}}L^2}{y^{\frac{n}{2}}}\ll\frac{y}{L^2}\cdot
N^{\frac{n\theta_n}{2}-\frac{1}{n}}
L^{\frac{n-1}{n}B_n+2-\frac{nA_n}{2}}\ll \frac{y}{L^2}.
\end{align*}
Thus, in all cases,
\[
\sup_{\alpha\in\mathfrak m}|T(\alpha,x,y)|
\ll\frac{y}{L^2}.
\]
Hence, by~\eqref{Est_int_primes},
\begin{align}\label{eq:Jminor-bound}
|J(\mathfrak m)|
\le
\sup_{\alpha\in\mathfrak m}|T(\alpha,x,y)|
\int\limits_{\mathfrak m}
&\left|
V(\alpha,N_1,2H)
V(\alpha,N_2,2H)
\right|\,d\alpha\ll
\frac{y}{L^2}\cdot\frac{H}{L}
\ll
\frac{H^2}{N^{\frac{n-1}{n}}L^3}.
\end{align}
From~\eqref{eq:circle-reduced}, the formula for $J(\mathfrak M_1)$,
and the estimates~\eqref{eq:JM2-bound}, \eqref{eq:Jminor-bound},
we obtain the asymptotic formula~\eqref{eq:main}.

\end{document}